\documentclass[11pt]{article}
\usepackage{amsmath, amsthm, amssymb}    
\usepackage{bridges}			
\usepackage{graphicx}			
\usepackage[colorlinks=true, urlcolor=blue, citecolor=black, linkcolor=black]{hyperref}  
\usepackage{subcaption}			

\begin{document}
\title{ Platonicons: The Platonic Solids Start Rolling}

\author{Katherine A. Seaton\textsuperscript{1} and David Hirsch\textsuperscript{2}\\
\vspace{10pt}\\
\textsuperscript{1}{Dept. Mathematics and Statistics, La Trobe University, Australia; k.seaton@latrobe.edu.au}\\
\textsuperscript{2}{Independent artist; playidea@gmail.com}} 

\date{}


\maketitle

\thispagestyle{empty}

\begin{abstract}

We describe the construction of a new family of developable rollers based on the Platonic solids. In this way kinetic sculptures may be realised, with the Platonic solids quite literally in their heart. We also describe the strong way in which the Platonicons circumscribe the Platonic solids.

\end{abstract}


\maketitle
\section*{A New Family of Platonic Developable Rollers}
Recently we announced the construction and some properties of a new family of solids, the polycons \cite{hs}, which generalise the sphericon \cite{DH, St}.  Over the course of their motion on a plane, rolling in an amusing manner, their entire surface makes contact with it. The polycons are based on regular polygons, and are named for the multiple pieces of identical cones which comprise their surface. In principle, there are infinitely many polycons, as many as there are regular polygons. The polycons are examples of what we termed developable rollers, defined fully in \cite{hs}.

In this paper we introduce another family of developable rollers discovered by David Hirsch in 2017, which are based on the fundamental Platonic solids, and which we call collectively the Platonicons. Though the number of Platonic solids is only five, there are more than five members of this family, though not infinitely many. In line with the naming of the underlying Platonic solid, we call these the tetrahedcons, cubicons (hexahedcons), octahedcons, dodecahedcons and icosahedcons. Figure 1 shows  photographs of some Platonicons, one based on each of the Platonic solids.

\begin{figure}[h!tbp]
	\centering
	\includegraphics[height=1.2in]{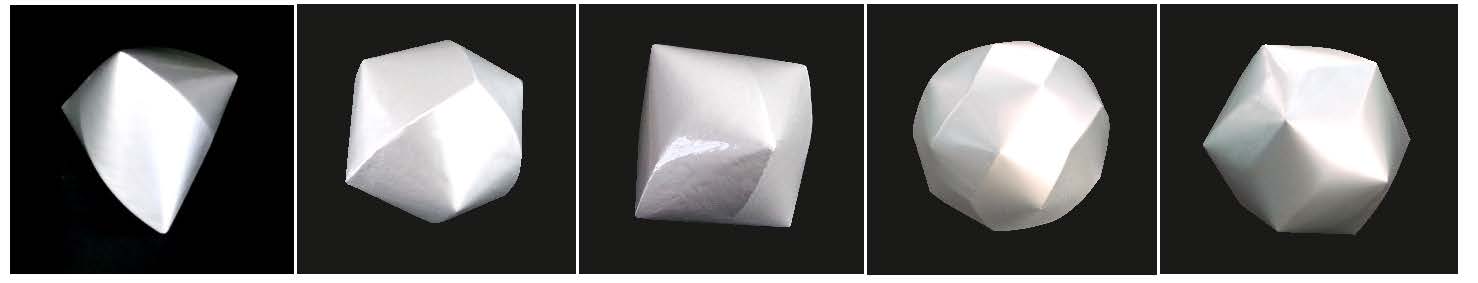}
	\caption{Five platonicons. 3D printed (plastic) and zinc paint. Indicative size: 10 $\times$ 8.5 $\times$ 10 cm.}
	
\end{figure}

\noindent From antiquity, the Platonic solids have enchanted artists and mathematicians. We hope that seeing them give rise to kinetic, rolling objects will inspire the artists and mathematicians among our readers.

\section*{Constructing the Platonicons}
\paragraph{Pieces of Cones}
Take any of the Platonic solids, with its vertices, edges and faces. The associated Platonicon is made up of the Platonic solid enhanced with identical modules made from  pieces of cones `glued' onto its faces. We distinguish one face of the solid and one vertex  of this face.  The modules are created thus:
\begin{enumerate}
\item[I] Construct a right circular cone with the distinguished vertex at its apex, and such that the edges connecting the adjacent vertices to the distinguished vertex are generators of the surface of the cone.  
\item[II] Construct a second identical cone, with its position determined by the nature of the distinguished face, as shown in Figure 2.

\begin{itemize}
 \item If the face is a triangle, this second cone has its apex at one of the vertices adjacent to the distinguished vertex. The edge connecting the vertices lies on both cones, and their curved surfaces cut each other in a conic section that lies above the altitude of the triangle that contains the third vertex. 
\item For the square face of the cube,  this second cone is placed so that its apex lies on the vertex diagonally opposite the distinguished vertex. The two cones cut each other in a conic section, lying above the \textit{other} diagonal of the square. 
\item Above the pentagonal face of the dodecahedron, not one but two further cones are constructed, with their apices at the two vertices which lie on the edge of the pentagon opposite to the distinguished vertex. For this case, there are three conic section edges associated with each module. 
\end{itemize}
In all cases, the volume common to the constructed cones and lying above the distinguished face comprises the module. A geometric argument based on the relative orientation of the solid's symmetry axes and the faces of its \textit{dual} Platonic solid shows that the vertex angle of the cones is the dihedral angle of the dual solid; this results in an object that will roll smoothly with its centre of mass maintaining constant height.
\item[III] Construct $n-1$ more identical modules, where $n$ is the number of faces of the platonic solid.
\end{enumerate}
\begin{figure}[h!tbp]
	\centering
	\includegraphics[width=4in]{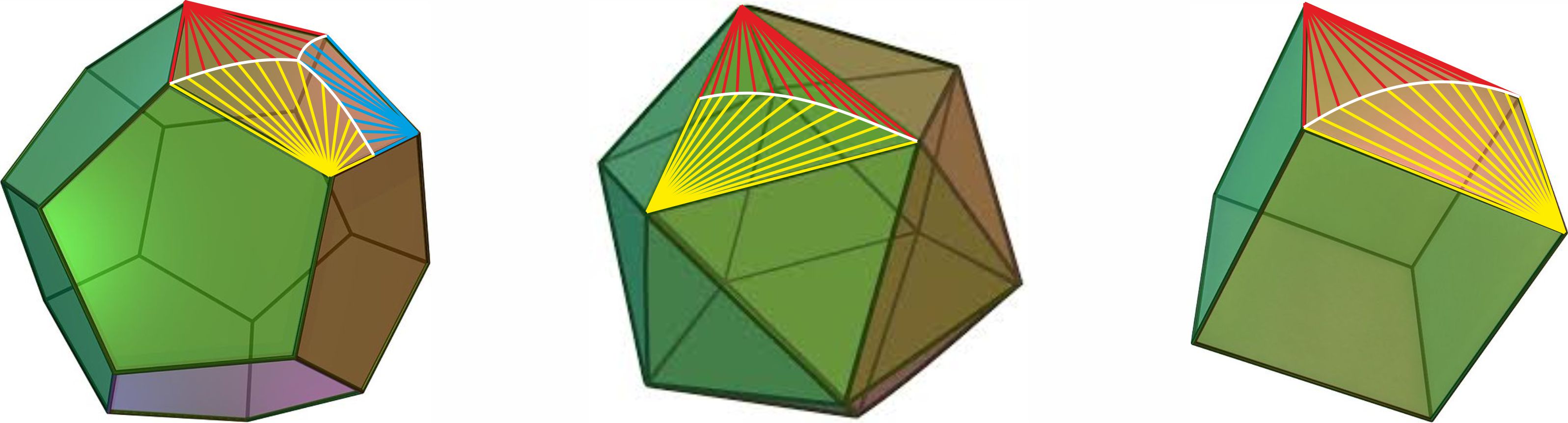}
	\caption{The modules constructed from pieces of cones, on a pentagonal, triangular and square face.}
	
\end{figure}

\paragraph{Orienting the Modules to Form the Platonicons}
Although there are 2, 3 or 5 orientations for each module (depending on the underlying polygon), not all of these result in distinct solids (up to symmetry and chirality). Significantly, not all are consistent with the conditions that give developable rollers. To create an object that can develop its entire surface as it rolls, the $n$ modules are constructed on the faces of the Platonic solid in such a way that the Platonicon has a single sinuous surface.  For example, in Figure 3 we indicate two cone arrangements on the tetrahedron which give Platonicons. 
When the path line shown in Figure 3 leaves the net, it must enter a different part of the surface but through the same edge and eventually create a closed path. Five octahedcons (three without an axis of rotational symmetry), three dodecahedcons and two of each of the other Platonicons have been identified; this has been done manually but systematically, marking up each of the Platonic solids. 
\begin{figure}[h!tbp]
	\centering
	\includegraphics[height=1.5in]{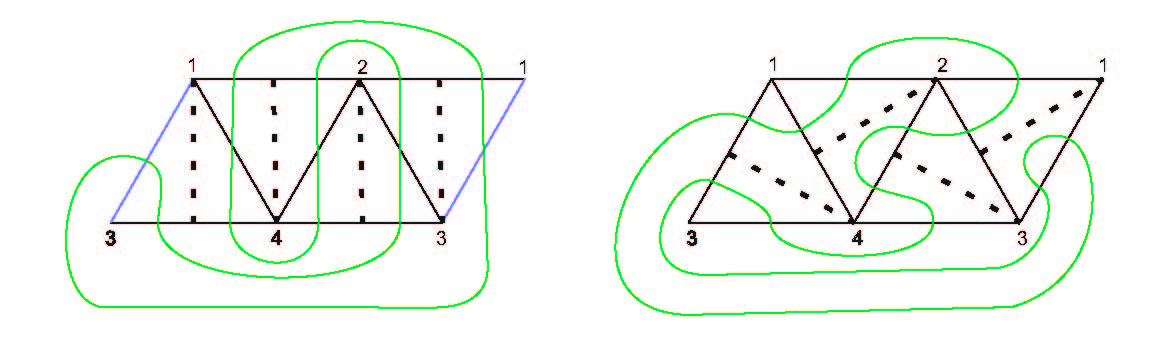}
	\caption{Tetrahedron nets showing features of the tetrahedcons. Dashed lines indicate the ridge where cone surfaces meet above the solid. Thin lines show the closed path of the rolling Platonicon. The vertices are numbered, with the same number indicating the same tetrahedcon vertex.}
	
\end{figure}
Animated images of some of these are available to view \cite{Think}. Figure 4 suggests why other orientations of the modules, or the construction of more cones, do not give rise to developable rollers. (We intend to expand on this in a future paper.)
\begin{figure}[h!tbp]
	\centering
	\includegraphics[height=1.5in]{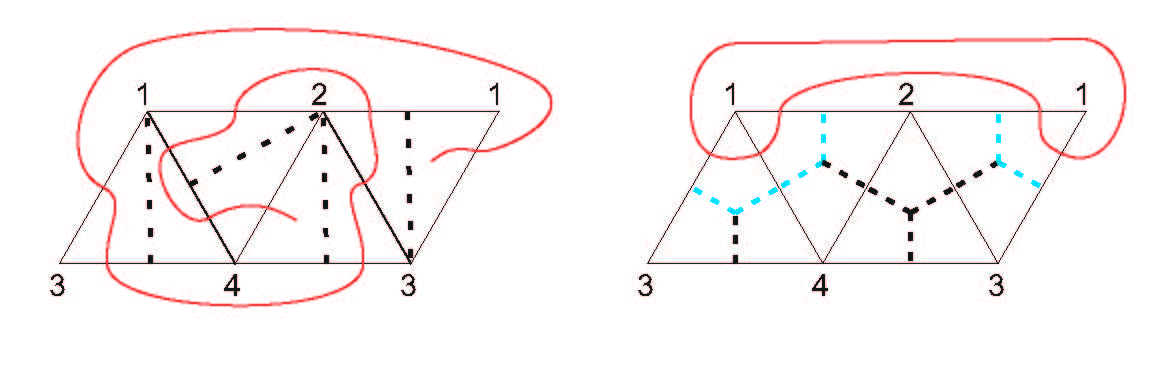}
	\caption{The module orientation on the left which does not allow a closed path, and modules made of three cones on the right which create closed regions of the surface, do not give  developable rollers.}
	
\end{figure}
\section*{Platonicon Properties}
\paragraph{Circumscribing the Platonic Solids}
Generally when one speaks of a Platonic solid being circumscribed, one thinks of the solid being inside a hollow sphere with only the vertices of the solid touching the sphere's surface, \textit {\`{a} la} Kepler's model of the planets. But a cylinder may also be circumscribed by a sphere: the two circular edges of the cylinder are latitudes (small circles) on the surface of the sphere. The Platonicons circumscribe the Platonic solids in \textit{both} ways. Not only are the vertices of the Platonic solid part of the Platonicon surface, but so too are its edges (see Figure 2). If a circumscribing sphere is imagined to kiss the Platonic solid, the Platonicon hugs it!
\paragraph{Vertices and Edges} A curved conic section edge of a Platonicon may end at a vertex which coincides with a vertex of the underlying Platonic solid, or it may be distinct,  and in some cases it merges into the surface of the Platonicon.  Two conic section edges may join each other to form a single edge, reducing the number of distinct edges. 
For generalised sphericons, Phillips \cite{AT2} investigated a connection with mazes. The edges of the Platonicons create even more complex mazes which warrant further investigation in future. 

\paragraph{Rolling Behaviour}
Like the other members of their extended family the developable rollers, the Platonicons roll smoothly but with meanderings in their motion, as first one and then another part of a cone's surface comes in contact with the surface upon which the rolling is taking place. 

However, some of the Platonicons, those described above as having a vertex that merges into their surface, exhibit  rolling behaviours not previously seen among the developable rollers.  By choosing the slope of the surface upon which they are rolling, they can be  encouraged to exhibit modes of rolling that do, or do not, develop their whole surface. Figure 5 shows four such modes for one of the tetrahedcons. To compare them in motion, which is perhaps more illuminating, in the videos of \cite{roll} see the cubicon which has only one mode of rolling, and the tetrahedcon or octahedcon which have several. 

\begin{figure}[h!tbp]
	\centering
	\includegraphics[height=1.1in]{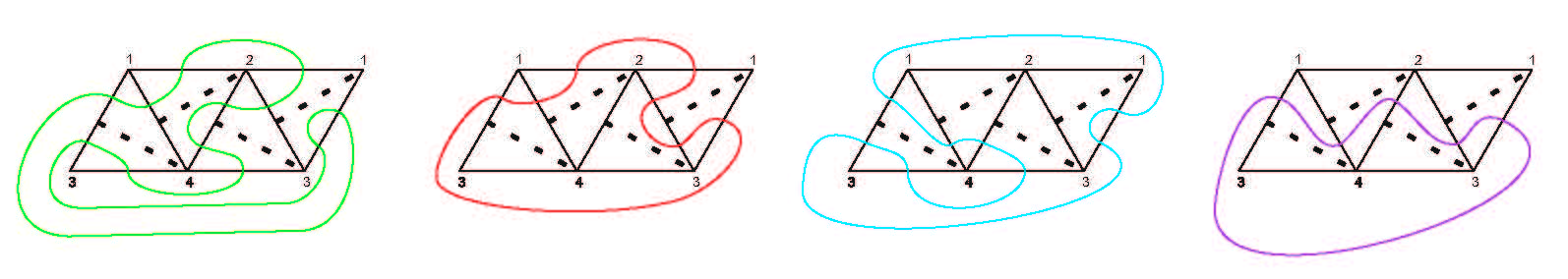}
	\caption{Four modes of rolling into which the same tetrahedcon can be encouraged. Only in the left-most of  these is the whole surface developed.}
	
\end{figure}

\section*{Summary and Conclusions}
We have introduced the Platonicons in this short paper, but we have only been able to touch on some of their features. We intend to continue to investigate these objects and their properties, and report on them. For instance, an obvious comparison to make will be that of the surface area and volume of each of the circumscribing Platonicons to the surface area and volume of the solid itself and that of the sphere that circumscribes the same Platonic solid in the customary (vertex only) way. The connection to mazes will be fascinating. The footprint of the Platonicons and the trajectory of their centers of mass should be studied. Constructions based on other solids, such as the cuboctahedron, and other Archimedean solids, will also be investigated. Recently, Daniel Walsh made us aware of his construction of `super advanced sphericons' \cite{dw}, which he also terms polyhedricons (various hexahedricons and the icosahedricon). There is potential to explore the similarities and differences between these and the Platonicons in consultation with him in future. 

As static objects, the Platonic solids have long inspired artists and mathematicians. We have shown in this paper that in the embrace of the Platonicons they can also dance.

\section*{Acknowledgements}
The underlying polyhedra in Figure 2 are based on those of https://en.wikipedia.org/wiki/User:Cyp.

    
{\setlength{\baselineskip}{13pt} 
\raggedright				

} 
   
\end{document}